\documentclass[a4paper,11pt]{article}
\usepackage[margin=25mm]{geometry}
\usepackage{amsmath,amssymb,amsthm,mathtools,bm}
\usepackage[hidelinks,unicode]{hyperref}
\usepackage{xurl,authblk,enumitem}
\setlist[enumerate]{leftmargin=2em,itemsep=0.2ex,topsep=0.4ex}

\hypersetup{pdftitle={An alternative proof of the SPR synthesis theorem for polynomial segments of arbitrary order},
  pdfsubject={Convex separation, contact certificates, and common-factor degree bounds}}
\numberwithin{equation}{section}
\theoremstyle{plain}
\newtheorem{theorem}{Theorem}[section]
\newtheorem{lemma}[theorem]{Lemma}
\newtheorem{proposition}[theorem]{Proposition}
\newtheorem{corollary}[theorem]{Corollary}
\theoremstyle{definition}
\newtheorem{definition}[theorem]{Definition}
\theoremstyle{remark}

\newcommand{\R}{\mathbb R}
\newcommand{\C}{\mathbb C}
\newcommand{\Repart}{\operatorname{Re}}
\newcommand{\conv}{\operatorname{conv}}
\newcommand{\lc}{\operatorname{lc}}
\newcommand{\K}{\mathcal K}
\title{An Alternative Proof of the SPR Synthesis Theorem\\
for Polynomial Segments of Arbitrary Order}

\author[1]{Zhengtai Xiong}
\author[1]{Tiancheng Xu}
\author[1]{Xiang Ji}
\author[1]{Guowei Dou}
\author[2,3]{Long Wang}
\author[1]{Wensheng Yu}
\affil[1]{School of Electronic Engineering, Beijing University of Posts and Telecommunications, Beijing 100876, China}
\affil[2]{Center for Systems and Control, College of Engineering, Peking University, Beijing 100871, China}
\affil[3]{Center for Multi-Agent Research, Institute for Artificial Intelligence, Peking University, Beijing 100871, China}
\date{}

\begin{document}
\maketitle

\begin{abstract}
For any positive integer $n$ and any two monic real Hurwitz polynomials of degree $n$, every convex combination of the two is Hurwitz stable if and only if there exists a common numerator of degree $n$ that makes both endpoint transfer functions strictly positive real. We give a unified proof based on convex separation, contact interpolation, and analytic properties of weakly positive real functions.
\end{abstract}

\section{Introduction}
In 1987, Dasgupta and Bhagwat \cite{dasgupta1987} studied robust strictly positive real (SPR) analysis and synthesis for interval polynomial families in the context of adaptive output-error identification. They asked how to choose a fixed numerator such that its ratio with every admissible denominator is SPR when the denominator coefficients vary within prescribed intervals. They showed that testing a given numerator against the entire family reduces to testing four distinguished denominator polynomials, and provided sufficient conditions and a construction for a common numerator. This problem connects design requirements in system identification with robust performance under coefficient uncertainty.

In 1990, Anderson et al.\ \cite{anderson1990} studied robust strict positive realness for $p/b$, where $p$ belongs to a convex polynomial family and $b$ is fixed, and obtained corresponding characterizations and construction conditions. A polynomial segment is generated by two endpoints. Its synthesis problem asks for one numerator that works for both endpoints and all their convex combinations. Hurwitz stability of the whole segment is necessary for such a numerator to exist; establishing sufficiency became a central question in segment synthesis.

In 2000, Wang and Yu \cite{wangyu2000} characterized SPR and weak SPR regions in numerator-coefficient space and proposed a robust synthesis method. Finding a common numerator could thus be viewed as finding an intersection of feasible regions associated with different denominators. This geometric approach laid the groundwork for subsequent results on polynomial segments of different orders.

In a 2002 preprint, Wang \cite{wang2002} gave a proof of the synthesis theorem for polynomial segments of arbitrary order. Yu and Wang \cite{yu2003} proved the sixth-order result in 2003. In the same year, Yu, Wang, and Ackermann \cite{yu2003ecc} presented the arbitrary-order result at the European Control Conference; their 2004 journal article \cite{yu2004} developed the SPR synthesis problem for polynomial families of arbitrary order. For two monic real Hurwitz polynomials of any positive degree, the entire segment is Hurwitz stable if and only if there is a common numerator of the same degree whose ratios with both endpoints are SPR. The same numerator then works for every polynomial on the segment. This theorem unifies the earlier results for individual low orders.

Building on the segment synthesis theorem in \cite{yu2003,yu2004}, we give another proof for arbitrary order. The geometric method of \cite{yu2004} obtains a candidate numerator of one lower degree from elliptic regions and then raises its degree by perturbation. Here we work directly in the coefficient space of numerators of the required degree. We express positivity at all frequencies as a compact convex-separation problem and extract a finite contact certificate at the first failure of synthesis.

Strict separation, Carath\'eodory's theorem, the maximum-modulus principle, and Lagrange interpolation are standard tools used here. The organizing step of our proof converts vanishing moments of the contact data into a lower bound on the degree of the common factor of the even and odd parts (Lemma~\ref{lem:residue}); weak positive realness gives an opposing upper bound. Simplicity of imaginary-axis zeros and the restriction on degree loss at infinity (Lemmas~\ref{lem:simple} and \ref{lem:drop}) make this comparison valid even when contact occurs at zero frequency or the leading coefficient of a supporting numerator vanishes. One interpolation lemma then excludes all four combinations of parity and leading-coefficient behavior, without separate tangent-plane equations for each degree. Strict separation finally yields a common numerator of degree $n$ that works at every parameter and frequency.

\section{Problem Statement and Proof Strategy}\label{sec:problem}

\subsection{Hurwitz stability and the main theorem}
\begin{definition}\label{def:hurwitz}
A nonzero polynomial $p\in\R[s]$ is \emph{Hurwitz stable} if
\begin{equation}
 p(z)=0\ \Longrightarrow\ \Repart z<0\qquad(z\in\C).\label{eq:hurwitz}
\end{equation}
For a positive integer $n$, let $H^n$ denote the set of all monic real Hurwitz polynomials of degree $n$. A polynomial family is \emph{robustly Hurwitz stable} if every member is Hurwitz stable.
\end{definition}

Fix a positive integer $n$, and write
\begin{align}
 a(s)&=s^n+a_1s^{n-1}+\cdots+a_n,\label{eq:a}\\
 b(s)&=s^n+b_1s^{n-1}+\cdots+b_n,\label{eq:b}\\
 p_\lambda(s)&=(1-\lambda)a(s)+\lambda b(s),\qquad0\leq\lambda\leq1,\label{eq:segment}\\
 \K&=\{p_\lambda:0\leq\lambda\leq1\}.\label{eq:K}
\end{align}
Set $a_0=b_0=1$. We adopt the same-degree SPR convention of \cite{yu2004}: for real polynomials $r,p$ both of degree $n$, the transfer function $r/p$ is \emph{strictly positive real} if $p$ is Hurwitz stable and
\begin{equation*}
 \Repart\frac{r(i\omega)}{p(i\omega)}>0\qquad(\omega\in\R).
\end{equation*}
Stability of $p$ ensures that this function has no pole in the closed right half-plane. We normalize the common numerator and denominators to be monic, so
\begin{equation*}
 \lim_{|\omega|\to\infty}\frac{r(i\omega)}{p(i\omega)}=1.
\end{equation*}
For sufficiently large $|\omega|$, the real part is at least $1/2$. On the remaining compact frequency interval, it is continuous and strictly positive, and therefore attains a positive minimum. Thus the real part has a uniform positive lower bound over all frequencies, and the frequency criterion also incorporates positivity at infinity.

\begin{theorem}\label{thm:main}
Let $a,b\in H^n$ be the monic real polynomials in \eqref{eq:a} and \eqref{eq:b}. There exists a monic real polynomial of degree $n$,
\begin{equation}
 \beta(s)=s^n+\beta_1s^{n-1}+\cdots+\beta_n,\label{eq:beta}
\end{equation}
such that both $\beta(s)/a(s)$ and $\beta(s)/b(s)$ are SPR if and only if
\begin{equation}
 (1-\lambda)a(s)+\lambda b(s)\in H^n
 \qquad(0\leq\lambda\leq1).\label{eq:main-iff}
\end{equation}
\end{theorem}

\begin{corollary}\label{cor:whole}
If $\K$ is robustly Hurwitz stable, the same numerator \eqref{eq:beta} can be chosen so that
\begin{equation}
 \Repart\frac{\beta(i\omega)}{p(i\omega)}>0
 \qquad(p\in\K,\ \omega\in\R).\label{eq:whole-spr}
\end{equation}
\end{corollary}

\subsection{Frequency-domain reduction}
For a polynomial $r(s)=\sum_{j=0}^nr_js^{n-j}$ written in descending powers, define
\begin{equation}
 E_r(t)=\sum_{j=0}^{\lfloor n/2\rfloor}(-1)^jr_{n-2j}t^j,\qquad
 O_r(t)=\sum_{j=0}^{\lfloor(n-1)/2\rfloor}(-1)^jr_{n-2j-1}t^j.
 \label{eq:evenodd}
\end{equation}
Then
\begin{equation}
 r(s)=E_r(-s^2)+sO_r(-s^2),\qquad
 r(i\omega)=E_r(\omega^2)+i\omega O_r(\omega^2).\label{eq:decomposition}
\end{equation}
Here $E_r(\omega^2)$ is the real part of $r(i\omega)$, while $\omega O_r(\omega^2)$ is its imaginary part. These arise from the even- and odd-power terms of $r$, respectively. Thus the even and odd parts of $r$ are $E_r(-s^2)$ and $sO_r(-s^2)$.

To search in coefficient space, for $v=(v_0,\ldots,v_n)\in\R^{n+1}$ set
\begin{equation}
 \beta_v(s)=v_0s^n+v_1s^{n-1}+\cdots+v_n.
\end{equation}
Here $v_j$ is the coefficient of $s^{n-j}$; the even and odd parts of $\beta_v$ are obtained from \eqref{eq:evenodd} by replacing $r_j$ with $v_j$. For any numerator $r$ of degree at most $n$, define
\begin{equation}
 Q_{p,r}(t)=E_p(t)E_r(t)+tO_p(t)O_r(t).\label{eq:cross}
\end{equation}
When $r=\beta_v$, abbreviate $Q_{p,\beta_v}$ to $Q_{p,v}$. This is the numerator of the real-part condition. If $p$ is Hurwitz stable, then
\begin{equation}
 \Repart\frac{\beta_v(i\omega)}{p(i\omega)}
 =\frac{Q_{p,v}(\omega^2)}{|p(i\omega)|^2}.\label{eq:realpart}
\end{equation}
Moreover, $\deg Q_{p,v}\leq n$ and $[t^n]Q_{p,v}=v_0$. If $v_0=0$, then $\deg Q_{p,v}\leq n-1$.

\begin{lemma}\label{lem:endpoints}
For fixed $v,t$,
\begin{equation}
 Q_{p_\lambda,v}(t)=(1-\lambda)Q_{a,v}(t)+\lambda Q_{b,v}(t).
 \label{eq:affine}
\end{equation}
If a numerator of degree $n$ satisfies strict frequency positivity for a monic stable denominator of degree $n$, its leading coefficient is strictly positive; division by that coefficient therefore makes it monic.
\end{lemma}
\begin{proof}
The first assertion follows from the affine dependence of $E_p,O_p$ on the denominator coefficients. For the second, let $c\ne0$ be the leading coefficient of the numerator. Then
\begin{equation*}
 \lim_{|\omega|\to\infty}\frac{\beta(i\omega)}{p(i\omega)}=c.
\end{equation*}
Positive real part at every frequency gives $c\geq0$, and exact degree $n$ gives $c\ne0$; hence $c>0$. Positive rescaling preserves the SPR condition.
\end{proof}

\subsection{Necessity and the strategy for sufficiency}
\begin{lemma}\label{lem:zeroexclusion}
Let $d_\tau$ be a family of monic real polynomials of degree $n$ whose coefficients depend continuously on $\tau\in[0,1]$. If $d_0$ is Hurwitz stable and no $d_\tau$ has an imaginary-axis zero, then every $d_\tau$ is Hurwitz stable.
\end{lemma}
\begin{proof}
The coefficients have a common bound $M$, so the Cauchy root bound gives $|z|\leq1+M$ for every root. Let $A$ be the set of parameters for which $d_\tau$ has a root in the open right half-plane. If $d_{\tau_0}(z_0)=0$ with $\Repart z_0>0$, choose $0<\varepsilon<\Repart z_0$. For $\tau$ near $\tau_0$, we have $|d_\tau(z_0)|<\varepsilon^n$. Factoring the monic polynomial into linear factors shows that some root of $d_\tau$ lies within $\varepsilon$ of $z_0$; otherwise the modulus of the product would be at least $\varepsilon^n$. Thus $A$ is relatively open.

If $\tau_j\in A$ converges to $\tau$, the uniform root bound allows us to select a convergent subsequence of right-half-plane roots. Its limit is a root of $d_\tau$ in the closed right half-plane. Since imaginary-axis zeros are excluded, its real part is still strictly positive. Hence $A$ is relatively closed. The interval is connected and $0\notin A$, so $A=\varnothing$. Excluding imaginary-axis zeros gives the claim. A fixed-degree family whose leading coefficients are continuous and strictly positive can first be normalized, yielding the same conclusion.
\end{proof}

The sufficiency proof has three steps.
\begin{enumerate}
 \item Compactify $t\in[0,\infty)$ to a closed interval and express existence of a common numerator as a strict-separation condition. If synthesis fails, obtain a nonzero supporting numerator $P$ and finite contact data at the first failure.
 \item Use weak positive realness of $P$ to control its imaginary-axis zeros and actual degree, giving an upper bound on the degree of the greatest common divisor of its even and odd parts.
 \item Merge contacts at equal frequencies. Coefficient balance yields vanishing moments, and Lagrange interpolation gives a lower bound on the same greatest-common-divisor degree. The bounds contradict each other, so strict separation holds.
\end{enumerate}
Sections~3--5 carry out these steps; Section~6 constructs the common numerator and proves the theorem.

\section{Frequency Compactification and a Critical Contact Certificate}\label{sec:certificate}
A critical contact certificate is positive-weight contact data arising at a failure of synthesis. It records nonnegative support, zero contacts, and coefficient balance.

\subsection{Compactified features and existence of a common numerator}
For the monomial direction $s^\ell$, set
\begin{equation}
 \phi_{p,\ell}(t)=
 \begin{cases}
  (-1)^ht^hE_p(t),&\ell=2h,\\
  (-1)^ht^{h+1}O_p(t),&\ell=2h+1.
 \end{cases}\label{eq:phi}
\end{equation}
Define the feature vector ordered by descending coefficient powers as
\begin{equation}
 \Phi_p(t)=(\phi_{p,n}(t),\phi_{p,n-1}(t),\ldots,\phi_{p,0}(t)).
\end{equation}
Then $v\cdot\Phi_p(t)=Q_{p,v}(t)$. Each component is a real polynomial of degree at most $n$. Put $u=t/(1+t)$ and compactify $q(t)=\sum_{j=0}^nq_jt^j$ by
\begin{equation}
 \widehat q(u)=\sum_{j=0}^nq_ju^j(1-u)^{n-j}.
\end{equation}
Write the corresponding vector as $\widehat\Phi_p(u)$. For $0\leq u<1$,
\begin{equation}
 \widehat\Phi_p(u)=(1-u)^n\Phi_p(u/(1-u)).
 \label{eq:compact}
\end{equation}
This polynomial vector is defined on the closed interval $[0,1]$. Let $\bm e_0,\bm e_n$ be the first and last standard basis vectors. Then
\begin{equation}
 \widehat\Phi_p(1)=\bm e_0,\qquad
 \widehat\Phi_p(0)=p_n\bm e_n.\label{eq:end-features}
\end{equation}
Consequently,
\begin{equation}
 v\cdot\widehat\Phi_p(u)>0\ (u\in[0,1])
 \quad\Longleftrightarrow\quad
 v_0>0\ \text{and}\ Q_{p,v}(t)>0\ (t\geq0).\label{eq:compact-pos}
\end{equation}

\begin{lemma}\label{lem:separation}
If $a,b$ are Hurwitz stable, let
\begin{equation}
 \mathcal C(a,b)=\conv(\widehat\Phi_a([0,1])\cup
                            \widehat\Phi_b([0,1])).
\end{equation}
A common monic numerator of degree $n$ that makes both ratios SPR exists if and only if $0\notin\mathcal C(a,b)$.
\end{lemma}
\begin{proof}
The union of the two compactified curves is compact, as is its convex hull in finite dimensions. If the origin is outside this convex hull, strict separation gives a vector $v$ whose inner product with every point of the hull is strictly positive. By \eqref{eq:compact-pos}, $v_0>0$ and both endpoint cross polynomials are strictly positive throughout $[0,\infty)$. The numerator $\beta_v/v_0$ then works. Conversely, the coefficient vector of a common monic numerator has strictly positive inner product with every point in the hull, so the hull cannot contain the origin.
\end{proof}

\subsection{First failure and finite contact data}
Suppose $\K$ is entirely stable but synthesis fails. Define
\begin{equation}
 \mathcal C_\lambda=\mathcal C(a,p_\lambda),\qquad
 D=\{\lambda\in[0,1]:0\in\mathcal C_\lambda\}.
\end{equation}
The set $D$ is nonempty because $1\in D$, whereas $0\notin D$ because $a/a=1$. If $\lambda_j\in D$ and $\lambda_j\to\lambda$, Carath\'eodory's theorem gives, for each $j$, a convex representation of the origin using at most $n+2$ generating points. Add zero weights to equalize the numbers of points. Frequencies and weights lie in compact sets, and there are only two endpoint labels. Passing to a common subsequence, the labels remain fixed while frequencies and weights converge. Continuity of the feature map in both parameter and frequency shows that the limiting combination still equals zero, hence $\lambda\in D$. Thus $D$ is closed, and the first failure parameter $\lambda_*=\min D$ lies in $(0,1]$.

\begin{lemma}\label{lem:certificate}
On the subsegment $[a,p_{\lambda_*}]$, there exist a nonzero supporting vector $v\in\R^{n+1}$, $m$ atoms with $m\leq n+1$, positive weights $w_j$, frequency parameters $u_j\in[0,1]$, and endpoint labels $d_j\in\{a,p_{\lambda_*}\}$ such that
\begin{align}
 &\sum_{j=1}^mw_j=1,\qquad
   \sum_{j=1}^mw_j\widehat\Phi_{d_j}(u_j)=0,\label{eq:atom-balance}\\
 &v\cdot\widehat\Phi_a(u)\geq0,\qquad
   v\cdot\widehat\Phi_{p_{\lambda_*}}(u)\geq0\quad(u\in[0,1]),\label{eq:support}\\
 &v\cdot\widehat\Phi_{d_j}(u_j)=0\quad(1\leq j\leq m).\label{eq:contact}
\end{align}
\end{lemma}
\begin{proof}
First, the origin lies on the boundary of $\mathcal C_{\lambda_*}$. If it were an interior point, choose sufficiently small positive and negative coordinate vectors in this convex hull and expand each as a finite convex combination of generating points. This yields points $x_1,\ldots,x_L$ on the generating curves whose convex hull contains a ball $B(0,\rho)$ for some $\rho>0$. Keep their frequencies and endpoint labels fixed, and denote the corresponding generating points at parameter $\lambda$ by $x_j(\lambda)$. For $\lambda<\lambda_*$ sufficiently close to $\lambda_*$, continuity gives $\max_j\|x_j(\lambda)-x_j\|<\rho/2$. For every unit vector $e$, $\max_j e\cdot x_j\geq\rho$, hence $\max_j e\cdot x_j(\lambda)>\rho/2$. If the origin were outside the convex hull of these perturbed points, strict separation would give a unit vector $e$ with $e\cdot x_j(\lambda)<0$ for all $j$, a contradiction. Therefore $0\in\conv\{x_j(\lambda)\}\subseteq\mathcal C_\lambda$, contrary to the definition of the first failure.

Boundary separation now gives a nonzero vector $v$ that is nonnegative on the hull. Every positive-weight atom in a convex representation of the origin must lie in the supporting hyperplane $v^\perp$. This hyperplane has dimension $n$, so Carath\'eodory's theorem there requires at most $n+1$ atoms. Removing any zero weights proves the claim.
\end{proof}

For the remainder of this section and for Sections~4 and 5, abbreviate $p_{\lambda_*}$ by $b$. The subsegment $[a,b]$ remains entirely stable; Section~6 returns to the original segment. Set $P=\beta_v$, $c=[s^n]P=v_0$, $F=E_P$, and $G=O_P$. Thus $P(s)=F(-s^2)+sG(-s^2)$, with even and odd parts $F(-s^2)$ and $sG(-s^2)$. The supporting condition gives $c\geq0$, and $Q_{a,P},Q_{b,P}$ are nonnegative on $[0,\infty)$.

\subsection{Merging equal-frequency contacts and moment equations}
Several contacts may occur at the same frequency. Interpolation requires distinct nodes, so we first merge equal-frequency contacts. The affine structure of the segment keeps each merged denominator stable. For $u_j<1$, put $t_j=u_j/(1-u_j)$ and merge equal finite $t_j$. Let $T\subset[0,\infty)$ be the set of distinct finite nodes. The merged positive weights $W_t$ and stable denominators are
\begin{equation}
 \begin{gathered}
 W_t=\sum_{j:t_j=t}w_j(1-u_j)^n>0,\qquad
 r_t=\frac{\displaystyle\sum_{j:t_j=t,\ d_j=b}w_j(1-u_j)^n}{W_t},\\
 d_t=(1-r_t)a+r_tb,\qquad 0\leq r_t\leq1.
 \end{gathered}\label{eq:grouped-weights}
\end{equation}
These sums include only finite-frequency atoms with $u_j<1$. All constituent weights are positive, so $r_t\in[0,1]$ and $d_t$ belongs to the stable segment. Affine dependence on the denominator coefficients gives
\begin{equation*}
 W_t\Phi_{d_t}(t)=\sum_{j:t_j=t}w_j(1-u_j)^n\Phi_{d_j}(t).
\end{equation*}
Let $\eta\geq0$ be the total weight of atoms at infinite frequency $u=1$. Then
\begin{equation}
 \sum_{t\in T}W_t\Phi_{d_t}(t)+\eta\bm e_0=0,\qquad
 \eta c=0.\label{eq:ordinary-balance}
\end{equation}
At each finite contact,
\begin{equation}
 F(t)E_{d_t}(t)+tG(t)O_{d_t}(t)=0\qquad(t\in T).\label{eq:group-contact}
\end{equation}
Define
\begin{equation}
 A_t=W_tE_{d_t}(t),\qquad B_t=W_t tO_{d_t}(t).\label{eq:AB}
\end{equation}
Taking all numerator-coefficient directions of degree less than $n$ in \eqref{eq:ordinary-balance} yields, for nonnegative integers $j$,
\begin{equation}
 \sum_{t\in T}t^jA_t=0\quad(2j<n),\qquad
 \sum_{t\in T}t^jB_t=0\quad(2j+1<n).\label{eq:moments}
\end{equation}
Infinite-frequency atoms contribute zero in these directions.

For $t>0$, the pair $(A_t,B_t)$ is nonzero because $d_t(i\sqrt t)\ne0$. At $t=0$, $A_0=W_0d_0(0)>0$ and $B_0=0$. There must be a positive finite node: if every atom were at $u=0$ or $u=1$, \eqref{eq:end-features} would make their positive-weight sum $\alpha\bm e_n+\eta\bm e_0$, with $\alpha,\eta\geq0$ not both zero. This cannot vanish because the two basis vectors are distinct for $n\geq1$. Write
\begin{equation}
 N=|T\cap(0,\infty)|\geq1,\qquad
 \nu=\begin{cases}1,&0\in T,\\0,&0\notin T,\end{cases}\qquad |T|=N+\nu.
 \label{eq:counts}
\end{equation}
Each node originated from a contact at an endpoint, so
\begin{equation}
 T\cap(0,\infty)\subseteq Z_a\cup Z_b,\qquad
 Z_d=\{t>0:Q_{d,P}(t)=0\}.\label{eq:cover}
\end{equation}
If $\nu=1$, the zero-frequency contact gives $P(0)=F(0)=0$.

\section{Zeros of a Weakly Positive Real Numerator and Common-Factor Degree}\label{sec:analytic}

\subsection{Extension to the half-plane and restrictions on zeros}
\begin{lemma}\label{lem:extension}
Let $d$ be a monic real Hurwitz polynomial of degree $n$, and let $P\in\R[s]$ have degree at most $n$ with $[s^n]P=c\geq0$. If
\begin{equation}
 \Repart\frac{P(i\omega)}{d(i\omega)}\geq0\qquad(\omega\in\R),
 \label{eq:boundaryweak}
\end{equation}
then
\begin{equation}
 \Repart\frac{P(z)}{d(z)}\geq0\qquad(\Repart z\geq0).
 \label{eq:closedweak}
\end{equation}
\end{lemma}
\begin{proof}
First, $P+\delta d$ is Hurwitz stable for every $\delta>0$. Consider the fixed-degree path $d+\tau P$, $0\leq\tau\leq1/\delta$. Its leading coefficient is $1+\tau c>0$, and
\begin{equation*}
 \Repart\frac{d(i\omega)+\tau P(i\omega)}{d(i\omega)}\geq1.
\end{equation*}
Thus no member of the path has an imaginary-axis zero. Since the initial polynomial $d$ is stable, Lemma~\ref{lem:zeroexclusion} and positive rescaling show that $P+\delta d$ is stable.

Apply the Cayley transform
\begin{equation}
 H(z)=\frac{P(z)-d(z)}{P(z)+d(z)}.\label{eq:cayley}
\end{equation}
The denominator $P+d$ is Hurwitz stable, so $H$ is analytic in the open right half-plane and continuous up to its boundary. From \eqref{eq:boundaryweak} and
\begin{equation*}
 |P+d|^2-|P-d|^2=4\Repart(P\overline d)
\end{equation*}
we obtain $|H(i\omega)|\leq1$. Comparison of leading terms gives $H(z)\to(c-1)/(c+1)$ as $z\to\infty$, whose modulus is at most $1$. Given $\varepsilon>0$, sufficiently large right semicircular arcs satisfy $|H|\leq1+\varepsilon$, while the imaginary-axis diameter satisfies $|H|\leq1$. The maximum-modulus principle on the right half-disk yields $|H(z)|\leq1+\varepsilon$. Letting $\varepsilon\downarrow0$ at a fixed point of the right half-plane gives $|H(z)|\leq1$. The same modulus identity then implies \eqref{eq:closedweak}.
\end{proof}

\begin{lemma}\label{lem:simple}
Under the hypotheses of Lemma~\ref{lem:extension}, if $P\ne0$, then $P$ has no zero in the open right half-plane, every imaginary-axis zero of $P$ is simple, and $P(x)\geq0$ for every $x\geq0$.
\end{lemma}
\begin{proof}
Put $f=P/d$. If $f$ had a zero $z_0$ in the open right half-plane, local factorization of a nonzero analytic function would give $f(z)=(z-z_0)^mh(z)$, with $m\geq1$ and $h(z_0)\ne0$. A direction near $z_0$ can be chosen along which the real part of this expression is negative, contradicting \eqref{eq:closedweak}. Hence there is no interior zero.

At an imaginary-axis zero $z_0=i\omega_0$, we have $d(z_0)\ne0$, so $f$ is analytic in a complex neighborhood of $z_0$. If the zero had multiplicity $m\geq2$, apply the same local factorization. Directions into the right half-plane have arguments in $(-\pi/2,\pi/2)$; multiplying this interval by $m$ allows a vector $\zeta$ with $\Repart\zeta>0$ and $\Repart(\zeta^mh(z_0))<0$. For sufficiently small $\rho>0$,
\begin{equation*}
 \Repart f(z_0+\rho\zeta)
 =\rho^m\Repart(\zeta^mh(z_0+\rho\zeta))<0,
\end{equation*}
again a contradiction. Thus every imaginary-axis zero is simple. Because $d$ has no such zero, $f$ and $P$ have the same multiplicity there.

Finally, a monic real Hurwitz polynomial satisfies $d(x)>0$ for $x\geq0$: its negative real roots contribute positive linear factors, and conjugate-root pairs contribute squared moduli. On the nonnegative real axis, \eqref{eq:closedweak} therefore yields $P(x)\geq0$.
\end{proof}

\begin{lemma}\label{lem:drop}
If the nonzero polynomial $P$ in Lemma~\ref{lem:simple} satisfies $[s^n]P=0$, then
\begin{equation}
 \deg P=n-1,\qquad\lc(P)>0.\label{eq:drop}
\end{equation}
\end{lemma}
\begin{proof}
Let $m=\deg P<n$ and form the fixed-degree reversals
\begin{equation}
 \widetilde P(w)=w^nP(1/w),\qquad
 \widetilde d(w)=w^nd(1/w).
\end{equation}
We have $\widetilde d(0)=1$, whereas $\widetilde P$ has a zero of multiplicity $n-m$ at $0$. For $w\ne0$, $\Repart w\geq0$ implies $\Repart(1/w)\geq0$, and $\widetilde P(w)/\widetilde d(w)=P(1/w)/d(1/w)$. Thus the reversed quotient is analytic near zero and has nonnegative real part on the right-half-plane side. The boundary-angle argument of Lemma~\ref{lem:simple} makes its zero at the origin simple. Hence $n-m=1$. Since $P$ is nonnegative and nonzero on the positive real axis, its leading asymptotic term also gives $\lc(P)>0$.
\end{proof}

\subsection{An upper bound on the common-factor degree}
A common root of the even and odd parts produces a symmetric pair of roots of $P$. The preceding zero properties therefore restrict common roots to distinct positive real numbers and make each such root a common contact of the two endpoints.

\begin{lemma}\label{lem:gcdroots}
Let $a,b\in H^n$ and $P\in\R[s]$ be nonzero with $\deg P\leq n$. Suppose $P$ has no zero in the open right half-plane and all its imaginary-axis zeros are simple. Write $F=E_P$, $G=O_P$, and let $J=\gcd(F,G)$ be monic, with $\gamma=\deg J$. Every complex root of $J$ is a distinct positive real root. If, in addition, $Q_{a,P}$ and $Q_{b,P}$ are nonnegative on $[0,\infty)$ and $Z_a,Z_b$ are their respective sets of distinct positive roots, then
\begin{equation}
 \gamma\leq |Z_a\cap Z_b|.\label{eq:gcd-overlap}
\end{equation}
\end{lemma}
\begin{proof}
If $J(r)=0$, then $F(r)=G(r)=0$. Choose $z\in\C$ with $z^2=-r$. By \eqref{eq:decomposition}, $P(z)=P(-z)=0$. Since $P$ has no open-right-half-plane zero, both $z$ and $-z$ must lie on the imaginary axis; hence $r\geq0$ is real. If $r=0$, then $F(0)=G(0)=0$, so $P(0)=P'(0)=0$, contrary to the simplicity of imaginary-axis zeros. Thus $r>0$.

If $J$ had a repeated root at $r>0$, then $(t-r)^2$ would divide both $F$ and $G$, and $(s^2+r)^2$ would divide $P(s)$, again contradicting simplicity. Hence the roots of $J$ are distinct positive real numbers, and their number equals $\deg J$. By Lemma~\ref{lem:cross-nonzero}, both endpoint cross polynomials are nonzero, so $Z_a$ and $Z_b$ are finite. Each root of $J$ is a root of both cross polynomials and thus belongs to $Z_a\cap Z_b$. This proves \eqref{eq:gcd-overlap}.
\end{proof}

\section{Vanishing Moments, Lagrange Interpolation, and Exclusion of Certificates}\label{sec:exclusion}

\subsection{Counting contacts of nonnegative cross polynomials}
\begin{lemma}\label{lem:cross-nonzero}
If $d$ is a monic real Hurwitz polynomial of degree $n$ and $P\in\R[s]$ is nonzero with $\deg P\leq n$, then $Q_{d,P}$ is not the zero polynomial.
\end{lemma}
\begin{proof}
The even--odd decomposition gives the polynomial identity
\begin{equation}
 P(s)d(-s)+P(-s)d(s)=2Q_{d,P}(-s^2).\label{eq:cross-odd-identity}
\end{equation}
If $Q_{d,P}\equiv0$, the left-hand side vanishes identically, so the rational function $P/d$ is odd. After cancellation, its poles can only lie in the open left half-plane. Oddness requires poles to occur in opposite pairs, so there are no finite poles. Since $P/d$ is proper, it must then be constant; oddness forces that constant to be zero, contrary to $P\ne0$.
\end{proof}

Every positive root of a nonzero polynomial that is nonnegative on $[0,\infty)$ has multiplicity at least two. Thus the number of distinct positive roots is at most half its degree; a root at zero uses at least one more degree. For the present certificate, let $C=|Z_a\cap Z_b|$ and $\gamma=\deg\gcd(F,G)$. From \eqref{eq:cover} and Lemma~\ref{lem:gcdroots},
\begin{equation}
 N+C\leq |Z_a|+|Z_b|,\qquad \gamma\leq C.\label{eq:overlap}
\end{equation}
Lemmas~\ref{lem:extension}--\ref{lem:gcdroots} apply because the supporting numerator is nonzero, its degree-$n$ coefficient is nonnegative, and its quotients with both stable endpoints are weakly positive real. Put
\begin{equation}
 d_* = \begin{cases}n,&c>0,\\n-1,&c=0.\end{cases}
\end{equation}
Both endpoint cross polynomials have degree at most $d_*$. Therefore
\begin{equation}
 N+C\leq2\lfloor\frac{d_*}{2}\rfloor.\label{eq:root-budget}
\end{equation}
Also $Q_{a,P}(0)=a_nP(0)$ and $Q_{b,P}(0)=b_nP(0)$. If $P(0)=0$, both cross polynomials vanish at zero, so
\begin{equation}
 N+C\leq2\lfloor\frac{d_*-1}{2}\rfloor.\label{eq:origin-budget}
\end{equation}

\subsection{A degree lower bound from vanishing moments}
The following lemma converts vanishing moments and contact conditions into a lower bound on the common-factor degree.
\begin{lemma}\label{lem:residue}
Let $M\geq1$, $k\in\mathbb N_0$, and let $t_1,\ldots,t_M$ be distinct real nodes. Let $z_i,y_i$ be real numbers such that the pairs $(z_i,y_i)$ are not all $(0,0)$, and let $E,O\in\R[t]$ satisfy $\deg E=k$ and $\deg O\leq k$. Choose $q\in\{0,1\}$ and suppose
\begin{align}
 &\sum_{i=1}^Mt_i^jz_i=0\quad(0\leq j<k),\label{eq:z-mom}\\
 &\sum_{i=1}^Mt_i^jy_i=0\quad(0\leq j<k+q),\label{eq:y-mom}\\
 &E(t_i)z_i+O(t_i)y_i=0\quad(1\leq i\leq M).\label{eq:res-contact}
\end{align}
Then
\begin{equation}
 2k+q<M+\deg\gcd(E,O).\label{eq:gap}
\end{equation}
\end{lemma}
\begin{proof}
Set $L(t)=\prod_{i=1}^M(t-t_i)$ and define the Lagrange interpolation polynomials
\begin{equation}
 U(t)=\sum_{i=1}^Mz_i\prod_{j\ne i}(t-t_j),\qquad
 V(t)=\sum_{i=1}^My_i\prod_{j\ne i}(t-t_j).\label{eq:residueUV}
\end{equation}
They satisfy $U(t_i)=z_iL'(t_i)$ and $V(t_i)=y_iL'(t_i)$. Distinctness of the nodes ensures $L'(t_i)\ne0$. Let $e_j$ be the $j$th elementary symmetric sum of $t_1,\ldots,t_M$, and let $e_j^{(i)}$ be the corresponding sum with $t_i$ omitted. For $0\leq j<M$,
\[
 [t^{M-1-j}]U=(-1)^j\sum_i z_i e_j^{(i)},
 \qquad
 e_j^{(i)}=\sum_{h=0}^j(-t_i)^h e_{j-h}.
\]
Thus the first $k$ vanishing moments make the top $k$ coefficients of $U$ vanish; similarly, the first $k+q$ vanishing moments eliminate the corresponding top coefficients of $V$. If $k\geq M$, the first $M$ equations in \eqref{eq:z-mom} form an invertible Vandermonde system and force all $z_i$ to vanish. Likewise, if $k+q\geq M$, the first $M$ equations in \eqref{eq:y-mom} force all $y_i$ to vanish. Consequently, in all cases,
\begin{equation}
 \deg U\leq M-k-1,\qquad \deg V\leq M-k-q-1.\label{eq:UVdegrees}
\end{equation}
We use the convention $\deg0=-\infty$. Since the interpolation data are not all zero, $U,V$ cannot both vanish. The polynomial $EU+OV$ has degree at most $M-1$ and, by \eqref{eq:res-contact}, vanishes at all $M$ distinct nodes. Hence
\begin{equation}
 EU+OV=0.\label{eq:EUOV}
\end{equation}
Write $J=\gcd(E,O)$, $E=JE_1$, and $O=JO_1$. Then $E_1,O_1$ are coprime and $E_1U+O_1V=0$, so $E_1\mid V$. Moreover, $V\ne0$: otherwise $EU=0$ would force $U=0$, a contradiction. With $\gamma=\deg J$, we obtain
\begin{equation}
 k-\gamma=\deg E_1\leq\deg V\leq M-k-q-1.
\end{equation}
Rearranging proves \eqref{eq:gap}. When $q=1$, the extra vanishing moment in the second group lowers the degree bound for $V$ by one more.
\end{proof}

Vanishing moments lower the degrees of the interpolation polynomials, while contact makes $E/\gcd(E,O)$ a divisor of $V$. The resulting lower bound on the common-factor degree can therefore be compared with the upper bound in \eqref{eq:overlap}.

\subsection{The four parity and leading-coefficient cases}
For $n=2k$, \eqref{eq:moments} gives the first $k$ vanishing moments of both $A$ and $B$. For $n=2k+1$, it gives the first $k+1$ moments of $A$ and the first $k$ moments of $B$. The contact condition is uniformly $F(t)A_t+G(t)B_t=0$. When $c=0$, Lemma~\ref{lem:drop} gives $\deg P=n-1$.

\begin{proposition}\label{prop:exclude}
For any $n\geq2$, a stable segment of monic polynomials of degree $n$ admits no critical certificate of the form \eqref{eq:atom-balance}--\eqref{eq:contact}.
\end{proposition}
\begin{proof}
Let $M=|T|=N+\nu$. By \eqref{eq:overlap}, $\gamma\leq C$ throughout. We consider four cases.

\emph{Case 1: $n=2k$ and $c>0$.}
Here $\deg F=k$ and $\deg G\leq k-1$. Apply Lemma~\ref{lem:residue} with $E=F$, $O=G$, $z=A$, $y=B$, and $q=0$ to obtain
\begin{equation}
 2k<N+\nu+\gamma.\label{eq:even-positive-gap}
\end{equation}
If $\nu=0$, \eqref{eq:root-budget} gives $N+C\leq2k$, an immediate contradiction. If $\nu=1$, zero-frequency contact gives $P(0)=0$. Equation~\eqref{eq:origin-budget} then gives $N+C\leq2k-2$, so $N+1+\gamma\leq2k-1$, again a contradiction.

\emph{Case 2: $n=2k$ and $c=0$.}
By \eqref{eq:drop}, $\deg P=2k-1$ and its leading coefficient is strictly positive. Hence $\deg G=k-1$ and $\deg F\leq k-1$. Interchange the even and odd roles in Lemma~\ref{lem:residue}: take $E=G$, $O=F$, $z=B$, $y=A$, replace $k$ by $k-1$, and set $q=1$. The first $k-1$ moments of $B$ and the first $k$ moments of $A$ vanish. Therefore
\begin{equation}
 2k-1<N+\nu+\gamma.\label{eq:even-zero-gap}
\end{equation}
After the degree drop, \eqref{eq:root-budget} gives $N+C\leq2k-2$. Since $\nu\leq1$, the right-hand side of \eqref{eq:even-zero-gap} is at most $2k-1$, a contradiction.

\emph{Case 3: $n=2k+1$ and $c>0$.}
Here $\deg G=k$ and $\deg F\leq k$. Set $E=G$, $O=F$, $z=B$, $y=A$, and $q=1$ to obtain
\begin{equation}
 2k+1<N+\nu+\gamma.\label{eq:odd-positive-gap}
\end{equation}
But \eqref{eq:root-budget} gives $N+C\leq2k$. Together with $\nu\leq1$, this makes the right-hand side at most $2k+1$, a contradiction.

\emph{Case 4: $n=2k+1$ and $c=0$.}
Here $n\geq3$, so $k\geq1$. By \eqref{eq:drop}, $\deg P=2k$, whence $\deg F=k$ and $\deg G\leq k-1$. Take $E=F$, $O=G$, $z=A$, $y=B$, and $q=0$, using only the first $k$ moments of each group. We obtain
\begin{equation}
 2k<N+\nu+\gamma.\label{eq:odd-zero-gap}
\end{equation}
If $\nu=0$, the post-drop bound $N+C\leq2k$ contradicts this inequality. If $\nu=1$, contact gives $P(0)=0$. Each endpoint cross polynomial, of degree at most $2k$, then has an additional zero at the origin. Hence $N+C\leq2k-2$ and the right-hand side is at most $2k-1$, again a contradiction.

These cases exhaust $c\geq0$, so a critical certificate cannot exist.
\end{proof}

\section{Proof of the Main Theorem}\label{sec:mainproof}
\begin{proof}[Proof of Theorem~\ref{thm:main}]
\textbf{Necessity.}
Suppose a monic real polynomial $\beta$ of degree $n$ makes both $\beta/a$ and $\beta/b$ SPR. Since $a,b\in H^n$, equation~\eqref{eq:realpart} gives $Q_{a,\beta}(t)>0$ and $Q_{b,\beta}(t)>0$ for every $t\geq0$. Thus \eqref{eq:affine} yields
\begin{equation}
 Q_{p_\lambda,\beta}(t)
 =(1-\lambda)Q_{a,\beta}(t)+\lambda Q_{b,\beta}(t)>0
 \quad(0\leq\lambda\leq1,\ t\geq0).\label{eq:necessity-cross}
\end{equation}
If some $p_\lambda$ had an imaginary-axis root $i\omega$, then $E_{p_\lambda}(\omega^2)=0$ and $\omega O_{p_\lambda}(\omega^2)=0$. By definition of the cross polynomial, this would give $Q_{p_\lambda,\beta}(\omega^2)=0$, contradicting \eqref{eq:necessity-cross}. Thus the entire segment has no imaginary-axis roots. All $p_\lambda$ are monic of degree $n$, and their coefficients depend continuously on $\lambda$. Lemma~\ref{lem:zeroexclusion}, applied from $p_0=a\in H^n$, gives $p_\lambda\in H^n$ for every $\lambda\in[0,1]$.

\textbf{Sufficiency.}
Suppose now that $p_\lambda\in H^n$ for every $\lambda\in[0,1]$. First let $n=1$. The stable denominators have the form $p_\lambda(s)=s+\alpha_\lambda$, where $\alpha_\lambda>0$. Choose any $h>0$ and set $\beta(s)=s+h$. Then
\begin{equation}
 \Repart\frac{i\omega+h}{i\omega+\alpha_\lambda}
 =\frac{\omega^2+h\alpha_\lambda}{\omega^2+\alpha_\lambda^2}>0
 \quad(0\leq\lambda\leq1,\ \omega\in\R).
\end{equation}
Now let $n\geq2$. If no common monic numerator of degree $n$ existed, Lemma~\ref{lem:separation} would give $0\in\mathcal C(a,b)$. Lemma~\ref{lem:certificate} would then produce a critical contact certificate on a stable subsegment of the original segment, contrary to Proposition~\ref{prop:exclude}. Hence $0\notin\mathcal C(a,b)$. Strict separation provides $v=(v_0,\ldots,v_n)\in\R^{n+1}$ such that $v\cdot\widehat\Phi_a(u)>0$ and $v\cdot\widehat\Phi_b(u)>0$ for every $u\in[0,1]$. Taking $u=1$ gives $v_0>0$. For $u<1$, \eqref{eq:compact-pos} gives $Q_{a,v}(t)>0$ and $Q_{b,v}(t)>0$ for all $t\geq0$. Define
\begin{equation}
 \beta(s)=\frac{\beta_v(s)}{v_0}
 =s^n+\frac{v_1}{v_0}s^{n-1}+\cdots+\frac{v_n}{v_0}.
 \label{eq:final-numerator}
\end{equation}
This is one fixed monic real polynomial of degree $n$. By affinity and \eqref{eq:realpart}, for every $\lambda\in[0,1]$ and $\omega\in\R$,
\begin{equation}
 \Repart\frac{\beta(i\omega)}{p_\lambda(i\omega)}
 =\frac{(1-\lambda)Q_{a,v}(\omega^2)+\lambda Q_{b,v}(\omega^2)}
        {v_0\,|p_\lambda(i\omega)|^2}>0.\label{eq:final-positive}
\end{equation}
Since $p_\lambda\in H^n$, the denominator does not vanish. Thus the same numerator makes the transfer function SPR along the whole segment. Taking $\lambda=0,1$ proves sufficiency and also establishes Corollary~\ref{cor:whole}.
\end{proof}

In summary, for every positive integer $n$ and any $a,b\in H^n$,
\begin{equation}
 \boxed{\begin{gathered}
 \exists\,\beta\in\R[s]\text{ monic with }\deg\beta=n:
 \quad \beta/a\text{ and }\beta/b\text{ are both SPR}\\
 \Longleftrightarrow\quad
 (1-\lambda)a+\lambda b\in H^n\quad(\forall\lambda\in[0,1]).
 \end{gathered}}\label{eq:complete-conclusion}
\end{equation}
The result has been formally verified in Lean.

\end{document}